\documentclass{amsart}
\usepackage{csquotes}
\usepackage[english]{babel}
\usepackage{url}
\usepackage{braket}
\usepackage{datetime}
\usepackage{graphicx}
\usepackage{float}
\usepackage[margin=0cm]{caption}
\usepackage{wrapfig}
\usepackage[font=small]{caption}
\usepackage{color}
\usepackage{subcaption}
\usepackage{upref}

\usepackage{mathtools}
\usepackage{amsmath}
\usepackage{amsfonts}
\usepackage{amssymb}
\usepackage{amsxtra}
\usepackage{amsthm}
\usepackage{enumerate}
\usepackage{mathrsfs}
\usepackage{dsfont}
\usepackage{bm}

\usepackage[backend=bibtex,style=alphabetic]{biblatex}
\usepackage{url}
\usepackage{pdfpages}

\usepackage[colorlinks=true,citecolor=blue,linkcolor=blue]{hyperref}

\newcommand{\Mod}[1]{\ \text{\normalfont mod}\ #1 }
\newcommand{\Z}{\mathbb{Z}}

\newcommand{\zn}{\Z_n}

\newcommand{\fp}{\mathbb{F}_p}
\newcommand{\zpm}{\Z_p^m}
\newcommand{\znm}{\Z_{n}^m}

\newcommand{\zpr}{\Z_{p^r}}
\newcommand{\zprm}{\Z_{p^r}^m}
\newcommand{\ideal}[1]{\langle #1 \rangle}
\newcommand{\field}{\zprm/\ideal{p}^m}

\newcommand{\bref}[2]{\hyperref[#2]{#1~\ref*{#2}}}
\newcommand{\brefdot}[2]{\hyperref[#2]{#1~\ref*{#2}}}
\newcommand{\bbref}[2]{\hyperref[#2]{#1~\ref*{#2}}}

\newcommand{\gf}{\mathscr{G}_f}

\newcommand\restr[2]{{
  \left.\kern-\nulldelimiterspace 
  #1 
  \vphantom{\big|} 
  \right|_{#2} 
  }}
  
\theoremstyle{plain}
\newtheorem{thm}{Theorem}[section]
\newtheorem{lem}[thm]{Lemma}
\newtheorem{proposition}[thm]{Proposition}
\newtheorem*{cor}{Corollary}

\newtheorem*{lemma-unnum}{}
\newtheorem*{prop-unnum}{}
\newtheorem*{namedthm}{Theorem \namedthmname}
\newcounter{namedthm}

\theoremstyle{definition}
\newtheorem{definition}{Definition}[section]
\newtheorem{examplex}{Example}[section]
\newenvironment{example}
  {\pushQED{\qed}\examplex}
  {\popQED\endexamplex}

\theoremstyle{remark}
\newtheorem{remarkx}{Remark}[section]
\newenvironment{remark}
  {\pushQED{\qed}\remarkx}
  {\popQED\endremarkx}

\makeatletter
\def\@settitle{\begin{center}%
  \baselineskip14\p@\relax
    \normalfont\bfseries\LARGE
\uppercasenonmath\@title
  \@title
  \end{center}%
}

\newenvironment{named}[1]
  {\def\namedthmname{#1}%
   \refstepcounter{namedthm}%
   \namedthm\def\@currentlabel{#1}}
  {\endnamedthm}

\renewcommand{\alpha}{e}
\renewcommand{\ell}{\mu}
\DeclareMathOperator{\im}{Im}

\makeatother

\title{Stabilization Bounds for Linear Finite Dynamical Systems}
\author{Bj\"orn Lindenberg}

\newdate{date}{11}{04}{2018}
\date{\displaydate{date}}

\begin{document}

\maketitle

\begin{abstract}
    \noindent A common problem to all applications of linear finite dynamical systems is analyzing the dynamics without enumerating every possible state transition. Of particular interest is the long term dynamical behaviour. In this paper, we study the number of iterations needed for a system to settle on a fixed set of elements. As our main result, we present two upper bounds on iterations needed, and each one may be readily applied to a fixed point system test. The bounds are based on submodule properties of iterated images and reduced systems modulo a prime. We also provide examples where our bounds are optimal. 
\end{abstract}

\section{Introduction}
\noindent A \emph{finite dynamical system} is an ordered pair $(X,f)$, which consists of a finite set~$X$ together with an iterating function $f$ over $X$. That is, elements of $X$ evolve by repeated application of $f$. Thus the study of a finite dynamical system is a study of sequences of the form 
\begin{equation}
\label{eq:intro1}
(  x,f(x),f^2(x),\cdots ),
\end{equation}
where $f$ is called the \emph{defining function}. In mathematical modelling, finite dynamical systems emerge in many different areas. Applications can be found in, for example, computational physics, electrical engineering, artificial intelligence, gene regulatory networks and cellular automata \cite{hernandez2005linear, colon2005boolean,colon2006monomial,bollman2007fixed}. Hence, $f$ often represents \emph{state transitions} or the \emph{passing of time}. \par
In this context, a \emph{linear finite dynamical system} $(M,f)$ is a system where $f$ is a linear map and $M$ is a \emph{finite $R$-module}. That is, $M$ is a generalized vector space with scalars in a finite commutative ring $R$, and is generated by a finite set of \emph{base elements}. Therefore, $f$ can be represented by a matrix $A$ over $R$, which describes the mapping of the base elements. \par 
An important problem in applications of finite dynamical systems, is the determination of long term behaviour without enumeration of all possible state transitions. The size of the underlying set $X$ may lead to iteration over all elements being computationally intractable. As an example, a finite dynamical system in the form of a Boolean modeling framework is demonstrated in \cite{laubenbacher2004computational}. There the framework models a gene regulatory network of 60 nodes, such that $|X| = 2^{60}$. Other frameworks for gene regulatory networks may have more than two states for each node, leading to a further increase in the size of the underlying set \cite{bollman2007fixed}. \par
When a linear finite dynamical system is given by $(\zpm, A)$, where $p$ is prime, the properties and long term behaviour is well understood \cite{bollman2007fixed,colon2005boolean,hernandez2005linear}. Although when dealing with more general linear systems $(\znm, A)$, one faces difficulties due to the lack of unique factorization \cite{xu2009linear,deng2015cycles}. 

\subsection{Fixed Point Systems}
Let $(X,f)$ be a finite dynamical system and consider an iterated sequence of the form given in \eqref{eq:intro1}. It is possible that after a certain number of iterations we reach a \emph{fixed point} $x_0$, such that $f(x_0)=x_0$. If every element in $X$ eventually iterates to a fixed point, not necessarily unique, then the system is called a \emph{fixed point system}. Otherwise elements of $X$ will eventually converge within a subset of $X$, which consists of fixed points and \emph{cycles}, where cycles are closed iteration loops of more than one element. An extensive study of cycles of linear finite dynamical systems can be found in \cite{deng2015cycles}. \par 
A \emph{fixed point system criterion} for $(\mathbb{F}_q^n, A)$, where $\mathbb{F}_q$ is the finite field with $q$ elements, is given by Bollman \emph{et al.} in \cite{bollman2007fixed}. It is based on the \emph{minimal polynomial} of $A$, i.e., the polynomial of least degree in $\mathbb{F}_q[x]$ which annihilates~$A$. As an alternative, G. Xu and Y.M. Zou presented another fixed point system criterion in \cite{xu2009linear}. Given $(R^m,A)$, where $R$ is a finite commutative ring, they showed that if the size of $R$ is a composite integer $n$, then the system is fixed point system if and only if
\begin{equation}
\label{eq:intro2}
  A^{k+1}=A^k,   
\end{equation}
where $k = \lceil m \log_2(n)\rceil$. Thus the integer $k$ is a derived bound on the number of iterations needed for the system to stabilize, and the criterion in \eqref{eq:intro2} provides an efficient way of determining if the system exhibits long term \emph{steady state} behaviour. \par
That every linear finite dynamical system eventually settles on fixed points and cycles is stated in \emph{Fitting's Lemma}, see \bbref{Theorem}{thm:fitting}. In the spirit of \cite{hernandez2005linear} we define the \emph{height} of a system as the minimum number of iterations needed for the system to settle. It is easy to show that the bound $k$ in \eqref{eq:intro2} can be replaced by any upper bound on the system height, generating an alternative fixed point system criterion. Thus, assuming a composite integer $n$ for the size of the ring, one could ask the question of how close $\lceil m \log_2(n)\rceil$ is to the actual system height, and if there is significant improvements to be found. \par
In this paper we give sharper height bounds. We also provide examples where these bounds are optimal.
Issues concerning implementation of our results and possible performance advantages is also addressed.\par
   We now proceed to discuss our main results in more detail along with examples, applications and analysis. 
   
\subsection{Main Results}
\noindent We present theorems which are height bounds for linear finite dynamical systems on free $\zn$-modules of the form $(\znm,f)$. As such they may be readily applied to the fixed point convergence test given in \eqref{eq:intro2}. The main results come in two varieties. One is independent of the defining function, and the other is based on the structure of smaller reduced systems. \par


\begin{named}{A}\text{\normalfont(A Function Independent Bound).}
\label{thm:a}
Let $(\Z_{n}^m, f)$ be a linear finite dynamical system, and let $s$ be the system height. If the prime factorization of $n$ is given by $p_1^{\alpha_1}p_2^{\alpha_2}\cdots p_{\omega}^{\alpha_{\omega}}$ then the height is bounded as
\begin{equation*}
   s \leq  m e_{\max},
\end{equation*}
where $e_{\max}$ is the largest of the prime factor exponents.
\end{named}

\noindent The function independent bound $m \alpha_{\max}$ of \bbref{Theorem}{thm:a}, where $\alpha_{\max}$ is the largest prime factor exponent of $n$, may be precalculated for any endomorphism given a certain module. Its proof relies on submodule properties of iterated images and reduced systems derived from the \emph{primary decomposition} of the underlying ring. 


\begin{named}{B}\text{\normalfont(A Function Dependent Bound).}
\label{thm:b}
Let $(\znm,f)$ be a finite dynamical system with a defining matrix $A$ over $\Z_n$, and let $s$ be the system height. If the prime factorization of $n$ is given by $p_1^{\alpha_1}p_2^{\alpha_2}\cdots p_{\omega}^{e_{\omega}}$, then the height is bounded as
\begin{equation*}
    s \leq  \max\left \{  s_1 \alpha_1,  s_2 \alpha_2, \dots, s_{\omega} \alpha_{\omega} \right \}, 
\end{equation*}
where each $s_i$ is the height of the reduced system $(\Z_{p_i}^m, A \Mod{p_i})$, $i = 1,2, \dots, \omega$.
\end{named}

\noindent The proof of \bbref{Theorem}{thm:b} is based upon properties derived from \emph{reduced} systems involving factor modules of prime powered orders. \par

\subsection{Examples and Analysis}
\begin{example}
With a system of the form $(\Z_{210}^{16}, f)$ we have $m = 16$ and $210 = 2^1 3^1 5^1 7 ^1$,
thus by \bbref{Theorem}{thm:a} the system height cannot be larger than $16\cdot1=16$. \par
Consider instead a system of the form $(\Z_{1960}^{4}, f)$. We have $m = 4$ and $1960 = 2^3 5^1 7^2$, yielding a maximal exponent of 3, as such the system height is less than or equal to $4\cdot3=12$.
\end{example}
\begin{remark}
\bbref{Theorem}{thm:a} gives an upper bound that is only dependent on the prime factor exponents. A system of the form $(\Z_6^3, f_1)$ has a height less than or equal to $3$, since $6 = 2^1 3^1$ and $3\cdot1=3$. The same bound applies to a system of the form $(\Z_{400827403}^3, f_2)$. Here we have the prime factorization
\begin{equation*}
   400827403 = 10333^1 38791^1,
\end{equation*}
thus in spite of $\Z_{400827403}^3$ being a set many orders of magnitude larger than $\Z_6^3$, the system will settle on cycles within three iterations.
\end{remark}
\begin{example}
Let $(\Z_{27720}^4,A)$ be a linear finite dynamical system, with a defining matrix
\begin{equation*}
      A = \begin{pmatrix}
 17453 & 19126 & 430 & 13601  \\
        3116 & 18264 & 19275 & 26452 \\
        22825 & 2401 & 22534 & 173 \\
        4496 & 13083 & 3885 & 12974
        \end{pmatrix},   
\end{equation*}
over $\Z_{27720}$. The prime factorization of 27720 is $2^3 3^2 5^1 7^1 11^1$. For each distinct prime $p$, we find the height of $(\Z_p^4, A \Mod{p})$ and multiply it with the corresponding exponent. This leads to a set of products
\begin{equation*}
  \left \{3\cdot3,2\cdot2,1\cdot0,1\cdot1,1\cdot0 \right \}=\left \{9,4,0,1,0 \right \}, 
\end{equation*}
and according to \bbref{Theorem}{thm:b}, the system will settle on cycles within 9 iterations. In this particular case the bound is actually optimal, i.e., one can show by numerical computations that 9 is the true height of the system. 
\end{example}
\begin{proposition}
Given an integer $n = p_1^{\alpha_1}p_2^{\alpha_2}\cdots p_{\omega}^{\alpha_{\omega}}$ and a finite dynamical system $(\znm,A)$. If $\ell_\text{A}, \ell_\text{B}$ are the calculated bounds from Theorem~\ref{thm:a} and \ref{thm:b}, then  
\begin{equation*}
    \ell_\text{B} \leq \ell_\text{A} \leq \lceil m \log_2(n) \rceil.
\end{equation*}
\end{proposition}
\begin{proof}
With a fixed prime $p_i$ in the factorization, we know from linear algebra that $\Z_{p_i}^m$ can have at most $m$ image reductions of dimension under the iteration of $A \Mod{p_i}$. Hence for each such system under the primary decomposition we have $s_i \leq m$, where $s_i$ is the height of the reduced system modulo $p_i$. This implies that
\begin{equation*}
    \ell_\text{B} = \max\left \{  s_1 \alpha_1,  s_2 \alpha_2, \dots, s_{\omega} \alpha_{\omega} \right \} \leq m \cdot \max\left \{  \alpha_1,  \alpha_2, \dots, \alpha_{\omega} \right \} = \ell_\text{A},
\end{equation*}
for the first inequality, while
\begin{align*}
&\ell_\text{A} = m \cdot  \max\left \{ \alpha_1, \dots, \alpha_{\omega}\right \}  \leq m \sum_{i=1}^{\omega} \alpha_i \leq m \sum_{i=1}^{\omega} \alpha_i \log_2 (p_i) \\
&= m \log_2(p_1^{\alpha_1}p_2^{\alpha_2}\cdots p_{\omega}^{\alpha_{\omega}}) = m \log_2(n) \leq \lceil m \log_2(n) \rceil
\end{align*}
shows the last inequality.
\end{proof}
\begin{figure}
\centering
\includegraphics[width=0.8\linewidth]{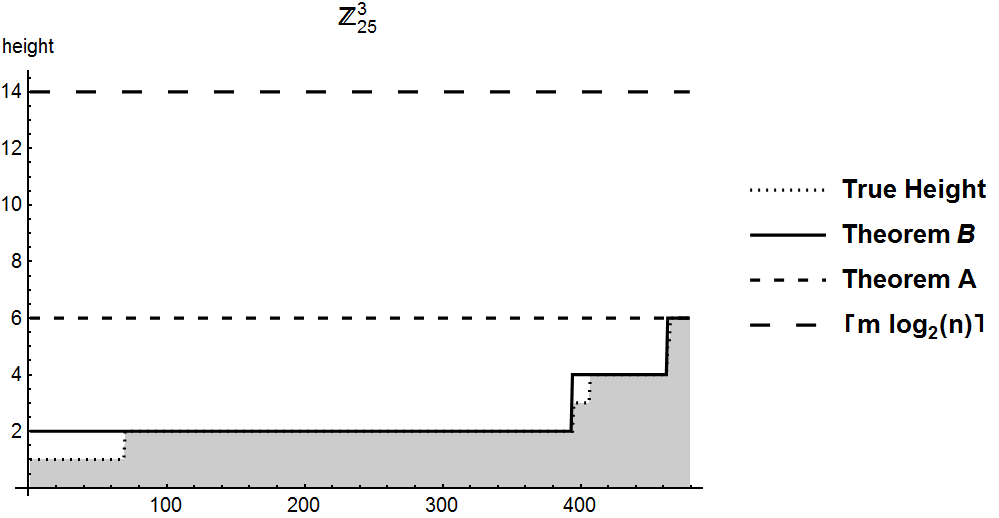}
\caption{A number of sampled linear systems $(\Z_{25}^3,f)$ for a bound comparison. Here \bbref{Theorem}{thm:a} demonstrates that it is a least upper bound over all singular matrices, and \bbref{Theorem}{thm:b} follows closely the actual height.}
\label{fig:h1}
\end{figure}
\begin{figure}
\includegraphics[width=0.8\linewidth]{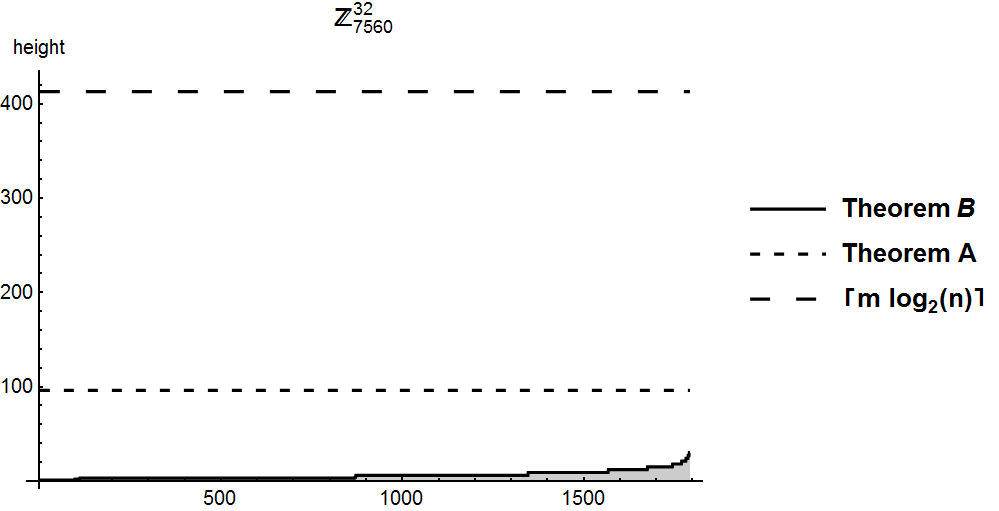}
\caption{A number of sampled linear systems $(\Z_{7560}^{32},f)$ for a bound comparison.}
\label{fig:h3}
\end{figure}
\begin{example}[Sampled Height Bounds] 
In \bbref{Figure}{fig:h1} height bounds are calculated for a number of sampled non-invertible matrices, and for reference the true height is included, which the data is sorted by. As can be seen, the bound of \bref{Theorem}{thm:b} follows closely the actual height of the systems, and \bref{Theorem}{thm:a} demonstrates that it gives a least upper bound concerning all endomorphisms. \bbref{Figure}{fig:h3} shows bounds for systems over a larger module. Here the number of distinct prime factors and base elements has increased compared to the former figure, which leads to relatively larger differences between the bounds. The figure shows that the function independent bounds may vastly overshoot the actual height of a system, where the height will be some integer less than or equal to
the function dependent bound of \bref{Theorem}{thm:b}. In both figures the function independent bound $\lceil m \log_2(n)\rceil$ derived in \cite{xu2009linear} is added for comparison.
\end{example}
\noindent Given a bound $k$ on the height, the efficient algorithm developed in \cite{xu2009linear} for the fixed point system test \eqref{eq:intro2} runs in time $O(\log(k) m^3)$, although here we assume naive matrix multiplication and faster algorithms for multiplication exists, see \cite{von2013modern}. When implementing \bbref{Theorem}{thm:b} in a fixed point convergence test there is an additional cost, since one needs to find the height $s_p$ of systems of the form $(\fp^m,A)$. One method is to use the minimal polynomial $m(x)$ of $A$, where $s_p$ is given by the first non-zero term in the expansion, such that $m(x) = x^{s_p} g(x)$, where $g(0) \neq 0$ \cite{hernandez2005linear}. This can be determined in run time $O(m^3)$ \cite{storjohann1998n}, giving us $O(\omega(n)m^3)$ bound computation complexity, where $\omega(n)$ is the number of distinct prime factors of $n$. Another method, which is inefficient but interesting in our case, is to use the fact that the reduced systems are defined over vector spaces. It follows that we can check the rank of successive powers of $A$ over $\fp$, such that $s_p$ is found when $\text{rank}(A^{s_p + 1}) =\text{rank}(A^{s_p})$. This will yield in total an $O(\omega(n) m^4)$ bound computation complexity.

In practice for smaller systems, using $\ell_\text{A}$ or $\lceil m \log_2(n)\rceil$ will be fast and efficient enough to outweigh any possible bound improvements gained by computing $\ell_\text{B}$. However given an application of a sufficiently large linear finite dynamical system the additional cost of computing the sharpest bound may be worth it. Given systems of the form $(\fp^m, A)$, we have
\begin{equation*}
   p^{t(m-t)}(p^m - 1)(p^m - p) \cdots (p^m - p^{t-1})
\end{equation*}
matrices that are bijective when restricted to their $t$-dimensional image subspace, such that $\text{rank}(A^2) = \text{rank}(A) = t$ \cite{laksov1994counting}. In particular, there are $\prod_{j=0}^{m-1}(p^m - p^j)$ invertible matrices, hence 
\begin{equation*}
 p^{m^2} - \prod_{j=0}^{m-1}(p^m - p^j) 
\end{equation*}
singular matrices in total. Therefore, the fraction of singular matrices that stabilizes in one iteration is given by
\begin{equation}
\label{eq:fraction}
\left ( \sum_{t = 0}^{m-1} p^{(m-t) t} Q_t \right ) / \left ( p^{m^2} - Q_m \right ),
\end{equation}
where $Q_0 = 1$, $Q_t = \prod_{j=0}^{t-1}(p^m - p^j)$. Since the numerator and denominator of~ \eqref{eq:fraction} are both monic polynomials in $p$ of degree $m^2 - 1$, the fraction approaches unity with increasing $p$. Furthermore, for a fixed prime, the ratio approaches some limiting value with increasing dimension, and for all except the smallest primes this limit corresponds to a considerable proportion of all singular matrices. To see this we fix a prime $p$ and let $q = 1/p \in (0,1/2]$. Next we introduce the $q$-Pochhammer symbol $(q;q)_n$ defined by
\begin{align*}
    (q;q)_0 &= 1, \\
    (q;q)_n &= \prod_{s =1}^{n} (1 - q^s),
\end{align*}
and we note that $Q_t$ can be expressed as
\begin{align*}
    Q_t &= \prod_{j = 0}^{t-1} (p^m - p^j) = q^{-m t} \frac{\prod_{s =1}^{m} (1 - q^s)}{\prod_{s =1}^{m-t} (1 - q^s)} = q^{- m t} \frac{(q;q)_m}{(q;q)_{m-t}}.
\end{align*}
Thus, for each term of ~\eqref{eq:fraction} we get
\begin{align*}
\frac{p^{m t - t^2} Q_t}{p^{m^2} - Q_m} &= \frac{q^{- 2 m t + t^2 + m^2}(q;q)_m}{(1-(q;q)_m)(q;q)_{m-t}} = \frac{(q;q)_m}{1-(q;q)_m}\frac{q^{(m-t)^2}}{(q;q)_{m-t}}, 
\end{align*}
which gives us the alternative form 
\begin{equation}
\label{eq:fraq_in_q}
    R_m = \frac{(q;q)_m}{1 - (q;q)_m} \sum_{s = 1}^{m} \frac{q^{s^2}}{(q;q)_s}
\end{equation}
for the ratio. \begin{figure}[t]
  \begin{subfigure}[t]{0.41\textwidth}
    \includegraphics[width=\textwidth]{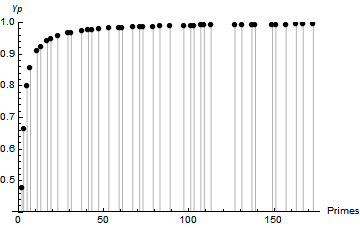}
    \caption{}
    \label{fig:fracp}
  \end{subfigure}
    \begin{subfigure}[t]{0.41\textwidth}
    \includegraphics[width=\textwidth]{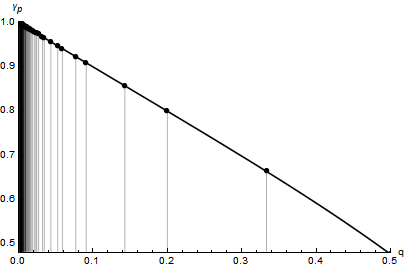}
    \caption{}
    \label{fig:fracq}
  \end{subfigure}
  \hfill
  \begin{subfigure}[t]{0.57\textwidth}
    \includegraphics[width=\textwidth]{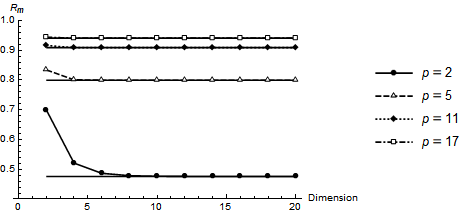}
    \caption{}
    \label{fig:fracm}
  \end{subfigure}
  \caption{The fraction $R_m$ of singular matrices over $\fp$ that stabilizes after one iteration. In (a) the value of a lower bound $\gamma_p$ over all dimensions is computed for increasing primes and in (b) for $q \in (0,1/2]$, which covers all primes with $q = 1/p$. In (c) the ratio is computed over fixed primes and variable dimensions along with computed lower bounds.}
\end{figure}
As a sequence, $R_m$ is strictly decreasing when $m \to \infty$, hence a good lower estimate is found by only counting the first three terms in the sum while sending $m$ to infinity. This yields a lower bound $\gamma_p$ over all dimensions such that
\begin{equation*}
    R_m > \gamma_p  = \frac{\phi(q)}{1 - \phi(q)}\left ( \frac{q}{(q;q)_1} + \frac{q^4}{(q;q)_ 2} + \frac{q^9}{(q;q)_3} \right ),
\end{equation*}
where $\phi(q) = (q;q)_\infty$ is the function used by Euler in the pentagonal number theorem \cite{ernst2012comprehensive}. Computed values for $\gamma_p$ on the first 40 primes is shown in \bbref{Figure}{fig:fracp}, and we can cover all primes by computing the values for $q \in (0,1/2]$, as shown in \bbref{Figure}{fig:fracq}. We can see that for all except the smallest primes, the matrices of height 1 make up a considerable proportion of the singular matrices regardless of dimension. In \bbref{Figure}{fig:fracm} computed ratios $R_m$ are shown for variable dimensions along with the lower bounds as solid lines.

Thus, for composite integers $n$ with relatively large prime factors, the probability of the true height being significantly less than $\ell_\text{A}$ is high for these systems. In fact, for systems of the form $(\zprm, A)$, where $A$ is singular and $p \gg 2$ is prime, the true height is most probably bounded by $1 \cdot r = r$ according to \bbref{Theorem}{thm:b} and the ratio $R_m$. Therefore, in such cases, any sufficient implementation may yield a performance advantage concerning fixed point system tests, since the actual computation of $\ell_\text{B}$ has a high probability of requiring only one or two multiplications with rank computations for each prime. Alternatively, given some large scale finite dynamical system, the reasoning above also shows that if we are using a function independent bound in the fixed point system test of \eqref{eq:intro2}, it is prudent to also test some intermediate values in the matrix power computation. This was suggested as a possible improvement in~\cite{xu2009linear}.  
\subsection{Related Work on Limit Cycles} When we have an upper bound $k$ on the height of $(\znm, f)$, we know by Fitting's Lemma that $N = f^k(\znm)$ is the submodule of $\znm$ where $f(N) = N$, i.e., these points make up the \emph{limit cycles} on which the dynamical system is stable and where $f$ acts like a bijection. An efficient algorithm to determine these cycles in soft-$O$ run time $m^3$ is presented by Wei~\textit{et~al}. in~\cite{wei2016dynamics}, where they use the bound and the automorphism $f_{\vert{N}}$. By the structure theorem of finitely generated abelian groups we can express $N$ by
\begin{equation}
\label{eq:relwork}
    N = \mathbb{Z}_{p_1^{a_1}} \oplus \mathbb{Z}_{p_2^{a_2}} \oplus \cdots \oplus \mathbb{Z}_{p_t^{a_t}},
\end{equation}
where the primes $p_i$ and positive integers $a_i$ are not necessarily distinct. Furthermore, the authors note that $N$ can be written as a direct sum of $f$-invariant submodules by grouping factors of \eqref{eq:relwork} having the same prime. This reduces the problem to a study of automorphisms on modules for each prime $p$ of the form
\begin{equation}
\label{eq:relwork2}
    \mathbb{Z}_{p^{a_1}} \oplus \mathbb{Z}_{p^{a_2}} \oplus \cdots \oplus \mathbb{Z}_{p^{a_r}},
\end{equation}
where $1 \leq a_1 \leq \cdots \leq a_r$. They show that such automorphisms in turn induces automorphisms on vector spaces, whose maximal cycle lengths can be determined with the aid of minimal polynomials, and these lengths will ultimately help us determine the possible cycle lengths of \eqref{eq:relwork2} under the corresponding automorphism. Solving for cycle elements of the known lengths for each invariant submodule of~\eqref{eq:relwork} then produces all possible cycles of $N$ by direct products. Note that a key ingredient here is knowing a bound $k$ on the true height such that $\im f^k = N$.
\subsection{Organization and Strategy}
We will follow the framework laid out in \cite{hernandez2005linear,bollman2007fixed,deng2015cycles} and start with some basic notions concerning vertices in the associated state space and structural properties of product systems in \bbref{\S}{sec:prelim}. \par
After the preliminaries, the proofs will be presented in \bbref{\S}{sec:proofa} and \bbref{\S}{sec:proofb}. First a small addition to the train of thought used by Xu and Zou in \cite[Theorem 2.1]{xu2009linear}. Iterated images of any linear map yields a descending chain of submodules, hence the minimum number of iterations is bounded by the maximum chain length over all functions. This produces an intermediary bound in \bbref{Lemma}{res11}, which involves $\Omega(n)$, the number of prime factors of $n$ counting multiplicity. Furthermore the main theorems are partially a consequence of \bbref{Lemma}{res12}, which equates the number of iterations needed for stabilization to that of subsystems arising from the primary decomposition of the underlying ring. \par

In \bbref{\S}{sec:proofb} we will derive properties of systems of the form $(\zprm, f)$. This is done by looking at smaller reduced systems $(\field, \overline{f})$ and $(\ideal{p}^m, f)$, where $\overline{f}$ is a function \emph{induced} by the natural homomorphism from the module to its factor module. This will lead us to a bound when $M = \zprm$ in \bbref{Lemma}{res16}, which combined with \bbref{Lemma}{res12} proves \bbref{Theorem}{thm:b}: A function dependent bound when $M = \znm$.

\section*{Acknowledgements}
\noindent I would like to thank the referee for valuable comments on how to improve this paper. I would also like to thank Karl-Olof Lindahl and Per-Anders Svensson at Linn\ae us University for productive discussions, pointers to literature, and for corrections and improvements of earlier drafts of this paper.
\noindent 

\section{Preliminaries}
\label{sec:prelim}
\noindent In a finite dynamical system $(X,f)$ every element is a vertex in the associated directed graph $\mathscr{G}_f$, called the \emph{state space}, whose edges are given by the ordered pairs $\set{(x,f(x)) | x \in X}$. We hereby define two different instances of vertices.

\begin{definition}\label{def:fds-cycle} 
Let $(X,f)$ be a finite dynamical system with an associated state space $\gf$. A \emph{cycle vertex} is an element $x_0$ in $X$, such that $f^l (x_0) =x_0$ for some positive integer $l$, i.e., there exists a path from $x_0$ to itself in $\gf$. The smallest such integer is called the \emph{period} of $x_0$. Furthermore, when cycle vertices are of period 1 they are called \emph{fixed points}.
\end{definition}

\begin{definition}\label{def:fds-leaf} 
Let $(X,f)$ be a finite dynamical system with an associated state space $\gf$. A \emph{leaf vertex} is an element $y$ in $X$, such that $y$ is a leaf in the associated state graph $\gf$, i.e., the equation $f(x)=y$ has no solution.
\end{definition}

\noindent Of special interest to us is the subgraph of $\gf$ consisting of only cycle vertices and their internal edges. When an element reaches a cycle vertex during iteration, the future iterated images will remain as cycle vertices. This is an immediate consequence of \bbref{Definition}{def:fds-cycle}, i.e., if $x_0$ is a cycle vertex then $f^k(x_0)$ is a cycle vertex of the same period for all nonnegative integers $k$.\par 
Thus given a finite dynamical system $(X,f)$ and an element $x$ in $X$, one could ask the questions: How many iterations does it take for $x$ to reach a cycle vertex in the associated state space? How many iterations does it take to have a guaranteed cycle vertex regardless of starting element? These quantities are defined as follows.

\begin{definition}
\label{def:fds-height}
Let $(X,f)$ be a finite dynamical system and let $x$ be an element in~$X$. The number of iterations it takes for $x$ to reach a cycle vertex during iteration is called the \emph{height} of $x$, and is denoted $h(x)$. 
\end{definition}

\begin{definition}
\label{def:fds-systemheight}
Let $(X,f)$ be a finite dynamical system. The \emph{system height} $s$, defined by
\begin{equation*}
    s = \max_{x \in X} \left \{ h(x) \right \},
\end{equation*}
is the largest possible height of all $x \in X$.
\end{definition}


\noindent When furnishing the set $\znm$ as a $\Z_n$-module and choosing any $m \times m$ matrix $A$ over $\Z_n$, we get a linear finite dynamical system $(\znm, f)$, where $f$ is a $\Z_n$-endomorphism defined by $A$. We now state a theorem which shows its \emph{primary decomposition} in the form of factor systems. The isomorphism ultimately has its root in the Chinese Remainder Theorem, and the primary decomposition of the underlying ring. 


\begin{thm}[\cite{deng2015cycles}]
\label{thm:lfdsiso}
Let $(\Z_{n}^m, A)$ be a linear finite dynamical system with a defining $m \times m$ matrix $A$ over $\Z_n$. If the prime factorization of $n$ is given by
\begin{equation*}
    n = p_1^{\alpha_1} p_2^{\alpha_2} \dots p_{\omega}^{\alpha_{\omega}},
\end{equation*}
then $(\Z_{n}^m, A)$ is isomorphic to
\begin{equation*}
   \Pi_{k=1}^{\omega}(\Z_{p_k^{\alpha_k}}^m, A \Mod{p_k^{\alpha_k}}), 
\end{equation*}
a product of linear finite dynamical systems. 
\end{thm}

\noindent Given a system of the form $(\zprm, f)$, where $p$ is prime, there is a natural reduction defined by the projective homomorphism to the quotient module $\zprm / \ideal{p}^m$. 

\begin{definition}
\label{def:indf}
Let $(\zprm, f)$ be a linear finite dynamical system, where $p$ is a prime number. Let $\ideal{p}^m$ be the submodule of $\zprm$ created by $m$ copies of the principal ideal $\ideal{p}$ in $\zpr$. Let 
\begin{equation*}
    \pi \colon \zprm \rightarrow \zprm / \ideal{p}^m \qquad x \mapsto x+\ideal{p}^m
\end{equation*}
be the natural homomorphism. We say that $\overline{f} \colon \zprm/\ideal{p}^m \rightarrow \zprm/\ideal{p}^m$ defined by $\overline{f} \pi = \pi f$, is the function \emph{induced by the natural homomorphism}, and where
\begin{equation*}
 \overline{f}(x+\ideal{p}^m) = f(x) + \ideal{p}^m   
\end{equation*}
for all $x+\ideal{p}^m \in \zprm/\ideal{p}^m$.
\end{definition}

\noindent By \emph{Fitting's Lemma} the state space of a linear system consists of copies of a single tree attached to cycle vertices. In particular, there will always be a copy of this tree in isolation converging to the zero element.

\begin{thm}[Fitting's Lemma \cite{bollman2007fixed}]
\label{thm:fitting}
Let $(M,f)$ be a linear finite dynamical system over commutative ring $R$. Then there exists a positive integer $s$ and submodules $N$ and $T$ such that
\begin{enumerate}[(i)]
    \item $N = f^s(M)$,
    \item $T = f^{-s}(0)$,
    \item $(M,f)$=$(N \oplus T,f|_N \oplus f|_T )$,
    \item $f|_N$ is invertible,
    \item $f|_T$ is nilpotent,
\end{enumerate}
where $f|_N$ and $f|_T$ are the restrictions of $f$ to each submodule.
\end{thm}

\begin{definition}
Let $(M,f)$ be a linear finite dynamical system over a commutative ring $R$. Let $T$ be the submodule of $M$, where for all $t \in T$, $f^k(t)=0$ for some nonnegative integer $k$. Then we will say that $T$ forms the \emph{nilpotent component} of the system.
\end{definition}


\section{Proof of the Function Independent Bound}
\label{sec:proofa}
\noindent In this section we will prove \bbref{Theorem}{thm:a}, a function independent bound on the height of a linear finite dynamical system of the form $(\znm, f)$. \par
\noindent Consider linear finite dynamical systems $(M, f)$ where $M = \znm$. We can derive an upper bound on the height of these systems which is independent of the defining function. By continuing the line of reasoning used in \cite{xu2009linear} we find the following lemma based on the prime factorization of $n$.


\begin{lem}
\label{res11}
Let $(\znm,f)$ be a linear finite dynamical system and denote by $s$ the height of the system. Then
\begin{equation*}
 s \leq m \Omega(n),
\end{equation*}
where $\Omega(n)$ is the number of prime factors of $n$ counting multiplicity.
\end{lem}

\begin{proof}
By the properties of the linear function every iterated image $f^k(\znm)$ must be a submodule of the previous iteration. Since a module is an abelian group under addition, each iterated image must also yield an additive subgroup of the previous iteration $f^{k-1}(\Z_n^m)$. By Lagrange's theorem $|f^{k}(\Z_n^m)|$ divides $|f^{k-1}(\Z_n^m)|$, i.e., the number of elements in the current iteration divides the number of elements in the previous iteration. \par
Let the prime factorization of $n$ be $p_1^{\alpha_1}p_2^{\alpha_2}\cdots p_{\omega}^{\alpha_{\omega}}$, where $\omega$ is the number of distinct prime factors of $n$. Then the size of the module can be expressed as $n^m=\sum_{i=1}^{\omega} p_i^{m \alpha_i}$. The longest possible chain of proper submodules will be for linear maps, whose iterations reduce in each step the number of elements with a factor of a single prime. This means that the maximum number of iterations $s$, for a chain of proper submodules is bounded above by
\begin{equation*}
s \leq \sum_{i=1}^{\omega} m \alpha_i = m \sum_{i=1}^{\omega} \alpha_i = m \Omega(n),
\end{equation*}
proving the lemma.
\end{proof}



\noindent \bbref{Theorem}{thm:lfdsiso} suggests that the height of a linear finite dynamical system $(\znm,f)$, is intimately connected to the possible heights of subsystems modulo $p^r$, which are derived from the primary decomposition. This is put on firm ground by the following lemma.

\begin{lem} 
\label{res12}
Let $(\Z_{n}^m, f)$ be a linear finite dynamical system with a defining matrix $A$ over $\Z_n$. If $n$ has the prime factorization $p_1^{\alpha_1}p_2^{\alpha_2}\cdots p_{\omega}^{\alpha_{\omega}}$, then the height of the system is given by the largest height of
\begin{equation*}
    \left \{ (\Z_{p_k^{\alpha_k}}^m,A \Mod{p_k^{\alpha_k}}) \right \},
\end{equation*}
the collection of subsystems derived from its primary decomposition.
\end{lem}

\begin{proof}
We want to show that for the product system given in \bbref{Theorem}{thm:lfdsiso}, the height is the largest height of the smaller factor systems. The result will then follow from the isomorphism. \par
The system is isomorphic to a product of subsystems of the form 
\begin{equation*}
   (\Z_{p_k^{\alpha_k}}^m,A \Mod{p_k^{\alpha_k}}), 
\end{equation*}
for each distinct prime factor of $n$. Let the product set be $X=X_1 \times X_2 \times \cdots \times X_\omega$, where $X_k=\Z_{p_k^{\alpha_k}}^m$. Let $\overline{f}$ be the function of the product system, such that
\begin{equation*}
    \overline{f}(x_1,x_2,\dots,x_\omega) = (f_1(x_1),f_2(x_2),\dots,f_\omega(x_\omega)),
\end{equation*}
and where each $f_k$ is defined by the matrix $A \Mod{p_k^{\alpha_k}}$. Let $s_1,s_2,\dots,s_\omega$ be the subsystem heights and let $s$ be the largest one. We have $s\geq s_k$ for each $k$ and according to Fitting's Lemma $f_k^{s+1}(X_k)=f_k^s(X_k)$ for all $k$. This is now the smallest possible nonnegative integer which settles each subsystem in the product system. Thus $s$ is the smallest nonnegative integer such that $\overline{f}^{s+1}(X)=\overline{f}^s(X)$, which yields $s$ as the system height of the product system. The result follows from the isomorphism.
\end{proof}

\noindent We can now state the proof of the function independent part of the main result. 


\begin{proof}[Proof of \normalfont \bbref{Theorem}{thm:a}]
By \brefdot{Lemma}{res12}, we know that the system height $s$ is given by 
\begin{equation*}
    s = \max \left \{s_1, s_2, \dots, s_\omega \right \},
\end{equation*}
where each $s_k$ is the height of a subsystem of the form $(\Z_{p_k^{\alpha_k}}^m,A \Mod{p_k^{\alpha_k}})$, derived from the primary decomposition. In addition, we know from \bref{Lemma}{res11}, that each $s_k$ is bounded above by $m \alpha_k$, where $\alpha_k$ is the prime factor exponent of $p_k$. Hence for each $s_k$, we have an upper bound $m \alpha_{\max}$, where $\alpha_{max}$ is the largest exponent in the prime factorization of $n$. Therefore we also have an upper bound on $s$ given by $m \alpha_{\max}$.
\end{proof}
\section{Proof of the Function Dependent Bound}
\label{sec:proofb}
\noindent In this section we will prove \bbref{Theorem}{thm:b}, a function dependent bound on the height of a linear finite dynamical system of the form $(\znm, f)$.\par
We will assume that the prime factorization of $n$ is known. By looking at factor subsystems, derived from the primary decomposition, and examining their links to corresponding reduced systems over fields, an upper bound is found for each subsystem. This bound will be function dependent and derived from the reduced system. The main result then follows when considering  \bbref{Lemma}{res12}, which states that the system height is given by the largest height of the factor systems. Hence by taking the largest derived bound of all subsystems, we get the result. \par
We now consider the relation between $(\zprm, f)$, where $p$ is prime, and a reduced system $(\field, \overline{f})$, where $\overline{f}$ is given in \bref{Definition}{def:indf}. Our ultimate goal is to derive from the reduction structural properties and height bounds for the larger system. We will start with a lemma showing a link between cycle vertices of both systems. This will provide a basis for other results and our understanding of how the state graphs are related.

\begin{lem}
\label{res2}
Let $(\zprm,f)$ be a linear finite dynamical system, and let $\ideal{p}$ be the maximal ideal in $\zpr$. Let $\overline{f}$ be the function induced by the natural homomorphism. Then $x_0+\ideal{p}^m$ is a cycle vertex in
\begin{equation*}
(\field, \overline{f})
\end{equation*}
if and only if the coset contains a cycle vertex in $(\zprm,f)$.
\end{lem}

\begin{proof}
We want to show that every cycle vertex $x+\ideal{p}^m$ in the reduced system, contains, as a coset, a cycle vertex in the larger system. We do this by demonstrating that for every element $x$ in the coset, there exists an element $w'$ in $\ideal{p}^m$, such that $x+w'$ is a cycle vertex in $(\zprm, f)$. Conversely, if the coset contains a cycle vertex $x_0$ such that $f^k(x_0)=x_0$ for some positive integer $k$, we demonstrate that the properties of the induced function must yield $x+\ideal{p}^m$ as a cycle vertex in $(\field, \overline{f})$. \par
$(\Rightarrow)$ Let $x+\ideal{p}^m$ be a cycle vertex in $(\field, \overline{f})$. Then there exists a $k\geq 1$ such that
\begin{equation*}
    (\overline{f})^k(x+\ideal{p}^m)=f^k(x)+\ideal{p}^m=x+\ideal{p}^m.
\end{equation*}
Hence $f^k(x)=x+w$ for some $w \in \ideal{p}^m$. Iterating in multiples of $k$ we get 
\begin{equation*}
f^{2k}(x)=f^k(f^k(x))=f^k(x+w)=f^k(x)+f^k(w)=x+w+f^k(w),
\end{equation*}
and in general for positive integers $q$ we have
\begin{equation*}
f^{qk}(x)=x+w+f^k(w)+\cdots+f^{(q-1)k}(w),
\end{equation*}
by the linearity of the function. We recall that $\ideal{p}^m$ is a submodule of $\zprm$, hence
\begin{equation*}
    w+f^k(w)+\cdots+f^{(q-1)k}(w)=w',
\end{equation*}
for some $w'$ in $\ideal{p}^m$. Choose an integer $q$ such that $qk>h(x)$, where $h(x)$ is the height of $x$. Then $f^{qk}(x) = x_0$, where $x_0$ is a cycle vertex in $(\zprm,f)$. We get
\begin{equation*}
    x_0=f^{qk}(x)=x+w',
\end{equation*}
which implies that $x_0$ is in $x + \ideal{p}^m$. \par
($\Leftarrow$) Assume that $x+\ideal{p}^m$ contains a cycle vertex $x_0$. Then there exists a $k\geq1$ such that $f^k(x_0)=x_0$. Hence
\begin{equation*}
    x+\ideal{p}^m=x_0+\ideal{p}^m = f^k(x_0)+\ideal{p}^m=(\overline{f})^k(x_0+\ideal{p}^m)=(\overline{f})^k(x+\ideal{p}^m)
\end{equation*}
is a cycle vertex in $(\field, \overline{f})$.
\end{proof}
\noindent A natural question is how the heights are related for elements of both systems. \bref{Lemma}{res2} gives us some information. An element $x$ cannot reach a cycle vertex~$y$ in $(\zprm, f)$ without turning $y+\ideal{p}^m$ into a cycle vertex in the reduced system. Hence $x$ cannot have a smaller height than $x+\ideal{p}^m$ in the reduced system. This effectively puts a lower bound on the height of all elements in the coset $x+\ideal{p}^m$. The following lemma makes it precise.


\begin{lem}
\label{res4}
Let $(\zprm,f)$ be a linear finite dynamical system, and let $\ideal{p}$ be the maximal ideal in $\zpr$. Let $\overline{f}$ be the function induced by the natural homomorphism. If $x+\ideal{p}^m$ has the height $t$ in 
\begin{equation*}
    (\field, \overline{f}),
\end{equation*}
then all $y \in x+\ideal{p}^m$ have heights larger than or equal to $t$ in $(\zprm, f)$.
\end{lem}

\begin{proof} We will do a proof by contradiction. Assume that there exists an element $y$ in $x+\ideal{p}^m$ of height $k$ less than $t$. Then $f^k(y) = y_0$, where $y_0$ is some cycle vertex in $(\zprm, f)$. Therefore, the $k$th iteration of $x+\ideal{p}^m$,
\begin{equation*}
    (\overline{f})^k(x+\ideal{p}^m)=f^k(x)+\ideal{p}^m = f^k(y)+\ideal{p}^m = y_0 + \ideal{p}^m,
\end{equation*}
contains a cycle vertex and according to \bref{Lemma}{res2}, the coset must be a cycle vertex. Hence $x+\ideal{p}^m$ iterates to a cycle vertex in less than $t$ steps. This contradicts the assumption.
\end{proof}

\noindent An immediate consequence of \bref{Lemma}{res4} is a lower bound on the height of the larger system by the height of the reduced system.

\begin{cor}
\label{cor1}
If $t$ is the height of $(\field, \overline{f})$, then the height $s$ of $(\zprm,f)$ is bounded below by $t$.
\end{cor}

\begin{proof}
If $t$ is the height of $(\field, \overline{f})$, then there exists an element $x+\ideal{p}^m$ with this height. By \bref{Lemma}{res4} every element of the coset must have a height equal to or larger than $t$. Hence for the maximal height $s$ of $(\zprm,f)$ we have $s\geq t$.  
\end{proof}

\noindent We now turn our attention to leaf elements. By Fitting's Lemma the system height can be found by the height of the nilpotent component, i.e., there exists a leaf element in the nilpotent component of maximum height. The next two lemmas show that we can deduce properties of $(\zprm, f)$ by examining leaf elements of the reduced system.\par 


\begin{lem}
\label{res5}
 Let $(\zprm,f)$ be a linear finite dynamical system, and let $\ideal{p}$ be the maximal ideal in $\zpr$. Let $\overline{f}$ be the function induced by the natural homomorphism. If $v+\ideal{p}^m$ is a leaf in
\begin{equation*}
    (\field, \overline{f}),
\end{equation*}
then the coset contains only leaf elements in $(\zprm, f)$.
\end{lem}

\begin{proof}
We will do a proof by contradiction.\par
Let $v+\ideal{p}^m$ be a leaf in $(\field, \overline{f})$. Assume that there exists a $y \in v+\ideal{p}^m$ and an $x \in \zprm$, such that $f(x) = y$. This implies that
\begin{equation*}
    \overline{f}(x+\ideal{p}^m)=f(x)+\ideal{p}^m=y+\ideal{p}^m=v+\ideal{p}^m,
\end{equation*}
contradicting the assumption that $v+\ideal{p}^m$ is a leaf.
\end{proof}

\begin{lem}
\label{res6}
 Let $(\zprm,f)$ be a linear finite dynamical system, and let $\ideal{p}$ be the maximal ideal in $\zpr$. Let $\overline{f}$ be the function induced by the natural homomorphism. If $v+\ideal{p}^m$ is a leaf in the nilpotent component of
\begin{equation*}
    (\field, \overline{f}),
\end{equation*}
then it contains a leaf in the nilpotent component of $(\zprm, f)$.
\end{lem}

\begin{proof}
Let $v+\ideal{p}^m$ be a leaf of height $t$ in the nilpotent component of $(\field, \overline{f})$. Then
\begin{equation*}
    (\overline{f})^t(v+\ideal{p}^m)=f^t(v)+\ideal{p}^m=\ideal{p}^m,
\end{equation*}
therefore $t$ is the smallest positive integer such that $f^t(v) = w$, for some $w \in \ideal{p}^m$. According to \bref{Lemma}{res5}, every $y \in v+\ideal{p}^m$ is a leaf, and by the corollary of \bbref{Lemma}{res4}, the height $s$ of $(\zprm,f)$ is bounded from below by $t$. \par 
Since $\ideal{p}^m$ is a submodule, $f^k(w)$ is an element in $\ideal{p}^m$ for all $w \in \ideal{p}^m$ and positive integers $k$. Therefore, given a system height of $s$, $f^s(v)=w_0$, where $w_0$ is some cycle vertex in $\ideal{p}^m$. \par 
If $w_0$ is of period $l$, i.e., $f^{l}(w_0)=w_0$, we choose an integer $q$ such that $r = q l - s$ is nonnegative. Then  $y = v-f^r(w_0) =v-w'$, is in $v+\ideal{p}^m$, and
\begin{equation*}
    f^s(y)=f^s(v-w')=f^s(v)-f^s(f^{ql-s}(w_0))=w_0-w_0=0,
\end{equation*}
proving the lemma.
\end{proof}

\noindent In dealing with the relevant reduced subsystems of $(\zprm, f)$, it is handy to work with isomorphic systems. The next proposition shows that the induced functions have simple representations based on the matrix of the linear function.


\begin{proposition} 
\label{res8}
Let $(\zprm,f)$ be a linear finite dynamical system with a defining matrix $A$ over $\zpr$. Let $\ideal{p}$ be the maximal ideal in $\zpr$ and let $\overline{f}$ be the function induced by the natural projective homomorphism. Then $ (\field, \overline{f})$ is isomorphic to $(\zpm, A \Mod{p})$, and $(\ideal{p}^m,f)$ is isomorphic to $(\Z_{p^{r-1}}^m, A \Mod{p^{r-1}})$.
\end{proposition}

\begin{proof} 
The proof is a simple examination of induced functions of bijective morphisms between systems. \par
Consider $(\field, \overline{f})$. Let $\pi \colon \zprm \rightarrow \field$ be the natural homomorphism, where 
\begin{equation*}
   \pi(x) =x+\ideal{p}^m = x \Mod{p} + \ideal{p}^m,
\end{equation*}
for all $x \in \zprm$. We may express the induced function $\overline{f}$ of $(\field, \overline{f})$ by
\begin{equation*}
    \overline{f}(x + \ideal{p}^m) =f(x)+ \ideal{p}^m = (Ax) \Mod{p} + \ideal{p}^m,
\end{equation*}
for all $x+\ideal{p}^m$. Consider the bijection
\begin{equation*}
   \psi \colon \field \rightarrow \zpm \qquad x +\ideal{p}^m \mapsto x \Mod{p},
\end{equation*}
which we take as a morphism between systems $(\field, \overline{f})$ and $(\zpm, \overline{g})$. Here we have an induced function
\begin{equation*}
   \overline{g}\colon \zpm \rightarrow \zpm \qquad \overline{g}(x)=\psi \overline{f} \psi^{-1} (x) = Ax \Mod{p},
\end{equation*}
for all $x \in \zpm$, and therefore $\psi$ acts as an isomorphism between $(\field, \overline{f})$ and $(\zpm, A \Mod{p})$. \par
Consider $(\ideal{p}^m,f)$. Here every element $w \in \ideal{p}^m$ can be represented uniquely by $p v$ for some element $v \in \Z_{p^{r-1}}^m$. Hence there exists a bijection $\phi: \ideal{p}^m \rightarrow \Z_{p^{r-1}}^m$, by 
\begin{equation*}
   \phi(p v) = v, 
\end{equation*}
which induces a function $g(v) = \phi f \phi^{-1}=\phi(A p v)=Av \Mod{p^{r-1}}$. As such, $(\ideal{p}^m,f)$ is isomorphic to $(\Z_{p^{r-1}}^m, A \Mod{p^{r-1}})$.
\end{proof}

\noindent The next lemma shows, that it is possible to put a lower and upper bound on the height of $(\zprm, A)$ by the height of the reduced system $(\zpm,A \Mod{p})$.

\begin{lem}
\label{res16}
Let $(\zprm,f)$ be a linear finite dynamical system with a defining matrix $A$ over $\zpr$. Let $s$ be the height of the system. Then 
\begin{equation*}
   s_1 \leq s \leq r s_1,
\end{equation*}
where $s_1$ is the height of the system $(\zpm, A \Mod{p})$.
\end{lem}

\begin{proof}
Let $z$ be a leaf of maximal height $s$ in the nilpotent component of $(\zprm,f)$. Since
\begin{equation*}
    f^s(z)+\ideal{p}^m = 0 + \ideal{p}^m = \ideal{p}^m,
\end{equation*}
$z+\ideal{p}^m$ is a coset in the nilpotent component of $(\field, \overline{f})$. As such $z$ reaches $\ideal{p}^m$ within $t$ iterations, where $t$ is the system height of the reduced system. Assuming maximal height, the reached element $w$ within $\ideal{p}^m$ is of system height $k$ of $(\ideal{p}^m,f)$. Hence $s$ is bounded by $t \leq s \leq t + k$. \par 
Working instead in the isomorphic systems of \brefdot{Proposition}{res8}, and denoting by $s_k$ the height of systems of the form $(\Z_{p^k}^m, A \Mod{p^k})$, we have $t = s_1$, $k = s_{r-1}$ and $s = s_r$. Hence by recursion
\begin{equation*}
    s_1 \leq s_r \leq s_1 + s_{r-1} \leq 2 s_1 + s_{r-2} \leq \cdots \leq r s_1,
\end{equation*}
proving the lemma.
\end{proof}

\noindent We are now ready to state the proof of the function dependent part of the main result. 

\begin{proof}[Proof of \normalfont \bbref{Theorem}{thm:b}]
As we know from \bref{Lemma}{res12}, the height is determined by the largest height of factor systems in the primary decomposition. Applying \bbref{Lemma}{res16} to each such factor system and taking the largest bound yields the result. \par 
Let $h_i, i=1,2,\dots,\omega$ be the heights of the factor subsystems $(\Z_{p_i^{\alpha_i}}^m,A \Mod{p_i^{\alpha_i}})$ derived from the primary decomposition. According to \bref{Lemma}{res16} each $h_i \leq \alpha_i s_i$, where $s_i$ is the height of $(\Z_{p_i}^m,A \Mod{p_i})$. Thus for all $i$, $h_i \leq \alpha_k s_k$, for some maximum product $\alpha_k s_k$ in 
\begin{equation*}
   \left \{ \alpha_1 s_1, \alpha_2 s_2, \dots, \alpha_{\omega} s_{\omega} \right \},
\end{equation*}
therefore we have a system height $s$ bounded by $a_k s_k$ according to \bbref{Lemma}{res12}.
\end{proof}



\section*{References}
\begingroup
\setlength{\emergencystretch}{2em}
\printbibliography[heading = none]
\endgroup



\end{document}